\documentclass[12pt]{article}

\usepackage{amsmath,amssymb,amsthm,amscd,a4wide,upref}

\usepackage[enableskew]{youngtab}

\usepackage{hyperref}
\hypersetup{colorlinks,citecolor=blue,filecolor=black,linkcolor=blue,urlcolor=blue}
\usepackage{enumerate}
\usepackage{mathtools}

\usepackage{lmodern}     
\usepackage[T1]{fontenc}

\DeclareSymbolFont{myletters}{OML}{ztmcm}{m}{it}
\DeclareMathSymbol{\uplambda}{\mathord}{myletters}{"15}

\begin{document}


\newcommand{\ad}{{\rm ad}}
\newcommand{\cri}{{\rm cri}}
\newcommand{\row}{{\rm row}}
\newcommand{\col}{{\rm col}}
\newcommand{\End}{{\rm{End}\ts}}
\newcommand{\Rep}{{\rm{Rep}\ts}}
\newcommand{\Hom}{{\rm{Hom}}}
\newcommand{\Mat}{{\rm{Mat}}}
\newcommand{\ch}{{\rm{ch}\ts}}
\newcommand{\chara}{{\rm{char}\ts}}
\newcommand{\diag}{{\rm diag}}
\newcommand{\non}{\nonumber}
\newcommand{\wt}{\widetilde}
\newcommand{\wh}{\widehat}
\newcommand{\ot}{\otimes}
\newcommand{\la}{\lambda}
\newcommand{\ls}{\ts\uplambda\ts}
\newcommand{\lam}{\uplambda}
\newcommand{\La}{\Lambda}
\newcommand{\De}{\Delta}
\newcommand{\al}{\alpha}
\newcommand{\be}{\beta}
\newcommand{\ga}{\gamma}
\newcommand{\Ga}{\Gamma}
\newcommand{\ep}{\epsilon}
\newcommand{\ka}{\kappa}
\newcommand{\vk}{\varkappa}
\newcommand{\si}{\sigma}
\newcommand{\vs}{\varsigma}
\newcommand{\vp}{\varphi}
\newcommand{\de}{\delta}
\newcommand{\ze}{\zeta}
\newcommand{\om}{\omega}
\newcommand{\Om}{\Omega}
\newcommand{\ee}{\epsilon^{}}
\newcommand{\su}{s^{}}
\newcommand{\hra}{\hookrightarrow}
\newcommand{\ve}{\varepsilon}
\newcommand{\ts}{\,}
\newcommand{\vac}{\mathbf{1}}
\newcommand{\di}{\partial}
\newcommand{\qin}{q^{-1}}
\newcommand{\tss}{\hspace{1pt}}
\newcommand{\Sr}{ {\rm S}}
\newcommand{\U}{ {\rm U}}
\newcommand{\BL}{ {\overline L}}
\newcommand{\BE}{ {\overline E}}
\newcommand{\BP}{ {\overline P}}
\newcommand{\AAb}{\mathbb{A}\tss}
\newcommand{\CC}{\mathbb{C}\tss}
\newcommand{\KK}{\mathbb{K}\tss}
\newcommand{\QQ}{\mathbb{Q}\tss}
\newcommand{\SSb}{\mathbb{S}\tss}
\newcommand{\TT}{\mathbb{T}\tss}
\newcommand{\ZZ}{\mathbb{Z}\tss}
\newcommand{\DY}{ {\rm DY}}
\newcommand{\X}{ {\rm X}}
\newcommand{\Y}{ {\rm Y}}
\newcommand{\Z}{{\rm Z}}
\newcommand{\Ac}{\mathcal{A}}
\newcommand{\Lc}{\mathcal{L}}
\newcommand{\Mc}{\mathcal{M}}
\newcommand{\Pc}{\mathcal{P}}
\newcommand{\Qc}{\mathcal{Q}}
\newcommand{\Rc}{\mathcal{R}}
\newcommand{\Sc}{\mathcal{S}}
\newcommand{\Tc}{\mathcal{T}}
\newcommand{\Bc}{\mathcal{B}}
\newcommand{\Ec}{\mathcal{E}}
\newcommand{\Fc}{\mathcal{F}}
\newcommand{\Gc}{\mathcal{G}}
\newcommand{\Hc}{\mathcal{H}}
\newcommand{\Uc}{\mathcal{U}}
\newcommand{\Vc}{\mathcal{V}}
\newcommand{\Wc}{\mathcal{W}}
\newcommand{\Yc}{\mathcal{Y}}
\newcommand{\Ar}{{\rm A}}
\newcommand{\Br}{{\rm B}}
\newcommand{\Ir}{{\rm I}}
\newcommand{\Fr}{{\rm F}}
\newcommand{\Jr}{{\rm J}}
\newcommand{\Or}{{\rm O}}
\newcommand{\GL}{{\rm GL}}
\newcommand{\Spr}{{\rm Sp}}
\newcommand{\Rr}{{\rm R}}
\newcommand{\Zr}{{\rm Z}}
\newcommand{\gl}{\mathfrak{gl}}
\newcommand{\middd}{{\rm mid}}
\newcommand{\ev}{{\rm ev}}
\newcommand{\Pf}{{\rm Pf}}
\newcommand{\Norm}{{\rm Norm\tss}}
\newcommand{\oa}{\mathfrak{o}}
\newcommand{\spa}{\mathfrak{sp}}
\newcommand{\osp}{\mathfrak{osp}}
\newcommand{\f}{\mathfrak{f}}
\newcommand{\g}{\mathfrak{g}}
\newcommand{\h}{\mathfrak h}
\newcommand{\n}{\mathfrak n}
\newcommand{\z}{\mathfrak{z}}
\newcommand{\Zgot}{\mathfrak{Z}}
\newcommand{\p}{\mathfrak{p}}
\newcommand{\sll}{\mathfrak{sl}}
\newcommand{\agot}{\mathfrak{a}}
\newcommand{\qdet}{ {\rm qdet}\ts}
\newcommand{\Ber}{ {\rm Ber}\ts}
\newcommand{\HC}{ {\mathcal HC}}
\newcommand{\cdet}{{\rm cdet}}
\newcommand{\rdet}{{\rm rdet}}
\newcommand{\tr}{ {\rm tr}}
\newcommand{\gr}{ {\rm gr}\ts}
\newcommand{\str}{ {\rm str}}
\newcommand{\loc}{{\rm loc}}
\newcommand{\Gr}{{\rm G}}
\newcommand{\sgn}{ {\rm sgn}\ts}
\newcommand{\sign}{{\rm sgn}}
\newcommand{\ba}{\bar{a}}
\newcommand{\bb}{\bar{b}}
\newcommand{\bi}{\bar{\imath}}
\newcommand{\bj}{\bar{\jmath}}
\newcommand{\bk}{\bar{k}}
\newcommand{\bl}{\bar{l}}
\newcommand{\hb}{\mathbf{h}}
\newcommand{\Sym}{\mathfrak S}
\newcommand{\fand}{\quad\text{and}\quad}
\newcommand{\Fand}{\qquad\text{and}\qquad}
\newcommand{\For}{\qquad\text{or}\qquad}
\newcommand{\OR}{\qquad\text{or}\qquad}
\newcommand{\grpr}{{\rm gr}^{\tss\prime}\ts}
\newcommand{\degpr}{{\rm deg}^{\tss\prime}\tss}

\numberwithin{equation}{section}

\renewcommand{\theequation}{\arabic{section}.\arabic{equation}}

\newtheorem{thm}{Theorem}[section]
\newtheorem{lem}[thm]{Lemma}
\newtheorem{prop}[thm]{Proposition}
\newtheorem{cor}[thm]{Corollary}
\newtheorem{conj}[thm]{Conjecture}
\newtheorem*{mthm}{Main Theorem}
\newtheorem*{mthma}{Theorem A}
\newtheorem*{mthmb}{Theorem B}
\newtheorem*{mthmc}{Theorem C}
\newtheorem*{mthmd}{Theorem D}

\theoremstyle{definition}
\newtheorem{defin}[thm]{Definition}

\theoremstyle{remark}
\newtheorem{remark}[thm]{Remark}
\newtheorem{example}[thm]{Example}
\newtheorem{examples}[thm]{Examples}

\newcommand{\bth}{\begin{thm}}
\renewcommand{\eth}{\end{thm}}
\newcommand{\bpr}{\begin{prop}}
\newcommand{\epr}{\end{prop}}
\newcommand{\ble}{\begin{lem}}
\newcommand{\ele}{\end{lem}}
\newcommand{\bco}{\begin{cor}}
\newcommand{\eco}{\end{cor}}
\newcommand{\bde}{\begin{defin}}
\newcommand{\ede}{\end{defin}}
\newcommand{\bex}{\begin{example}}
\newcommand{\eex}{\end{example}}
\newcommand{\bes}{\begin{examples}}
\newcommand{\ees}{\end{examples}}
\newcommand{\bre}{\begin{remark}}
\newcommand{\ere}{\end{remark}}
\newcommand{\bcj}{\begin{conj}}
\newcommand{\ecj}{\end{conj}}

\newcommand{\bal}{\begin{aligned}}
\newcommand{\eal}{\end{aligned}}
\newcommand{\beq}{\begin{equation}}
\newcommand{\eeq}{\end{equation}}
\newcommand{\ben}{\begin{equation*}}
\newcommand{\een}{\end{equation*}}

\newcommand{\bpf}{\begin{proof}}
\newcommand{\epf}{\end{proof}}

\def\beql#1{\begin{equation}\label{#1}}


\newcommand{\Res}{\mathop{\mathrm{Res}}}

\title{\Large\bf Classical $\Wc$-algebras for centralizers}

\author{{A. I. Molev\quad and\quad E. Ragoucy}}

\date{} 
\maketitle

\vspace{35 mm}

\begin{abstract}
We introduce a new family of Poisson vertex algebras $\Wc(\agot)$ analogous to the
classical $\Wc$-algebras. The algebra $\Wc(\agot)$ is associated with
the centralizer $\agot$ of an arbitrary nilpotent element in $\gl_N$.
We show that $\Wc(\agot)$ is an algebra of polynomials in infinitely many variables
and produce its free generators in an explicit form. This implies that
$\Wc(\agot)$ is isomorphic to the center at the critical level
of the affine vertex algebra associated with $\agot$.

\medskip

Report LAPTH-050/19

%

\end{abstract}


\vspace{45 mm}

\noindent
School of Mathematics and Statistics\newline
University of Sydney,
NSW 2006, Australia\newline
alexander.molev@sydney.edu.au

\vspace{7 mm}

\noindent
Laboratoire de Physique Th\'{e}orique LAPTh\newline
Universit\'{e} Grenoble Alpes, Universit\'{e} Savoie Mont Blanc, CNRS\newline
F-74000, Annecy, France\newline
eric.ragoucy@lapth.cnrs.fr

%

\newpage

\section{Introduction}
\label{sec:int}
\setcounter{equation}{0}

The {\em classical $\Wc$-algebras} associated with simple Lie algebras $\g$
were introduced by Drinfeld and Sokolov~\cite{ds:la} as
Poisson algebras of functions on infinite-dimensional manifolds.
A detailed review of the
constructions of \cite{ds:la} and
additional background of the theory can be found
in a more recent work by De Sole, Kac and Valeri~\cite{dskv:cw}.
As shown in \cite{dskv:cw},
the classical $\Wc$-algebras can be understood as
Poisson vertex algebras and can be used to
produce integrable hierarchies of bi-Hamiltonian equations.

The same framework of Poisson vertex algebras was used
in our paper \cite{mr:cw} to construct explicit generators
of the classical $\Wc$-algebras $\Wc(\g)$ associated to principal nilpotent elements
of simple Lie algebras $\g$
in types $A,B,C,D$ and $G$.
The generators are given in a uniform way as certain determinants
of matrices formed by elements of differential algebras associated with $\g$.
Another approach involving generalized quasideterminants
was developed in \cite{dskv:cwgln} to
describe generators of the classical $\Wc$-algebras
in type $A$ associated with arbitrary nilpotent elements of $\g$;
see also \cite{dskv:ca} for extensions of this approach to
other classical Lie algebras.

The principal classical $\Wc$-algebras also emerge from a different
perspective. Due to a theorem of Feigin and Frenkel~\cite{ff:ak},
\cite[Ch.~4]{f:lc}, the Poisson algebra $\Wc({}^L\g)$ associated to the Langlands dual Lie algebra ${}^L\g$
is isomorphic to the center at the critical level
of the affine vertex algebra $V(\g)$;
see also \cite[Ch.~13]{m:so} for a correspondence
of generators in the classical types.
In a recent work of Arakawa and Premet~\cite{ap:qm} a description of the
center at the critical level of the affine vertex algebra $V(\agot)$
was given for a family of centralizers $\agot=\g^e$ of certain
nilpotent elements $e\in\g$. Explicit generators of the center in type $A$
were constructed in \cite{m:cc}.

Our goal in this paper is to introduce analogues of the classical $\Wc$-algebras
for the centralizers $\agot=\g^e$ in the case $\g=\gl_N$.
Our classical $\Wc$-algebras $\Wc(\agot)$ turn out to fit the
general framework of Poisson vertex algebras of \cite{dskv:cw}
by possessing a $\lam$-bracket. We prove that $\Wc(\agot)$
is an algebra of polynomials in infinitely many variables and produce
a family of its free generators in an explicit form.
Moreover, we show that a Miura-type map leads to a natural isomorphism between
the algebra $\Wc(\agot)$
and the center of the affine vertex algebra $V(\agot)$ at the critical level.

\section{Definition of $\Wc(\agot)$}
\label{sec:def}

Suppose that $e\in\g=\gl_N$ is a nilpotent matrix with Jordan blocks of sizes
$\la_1,\dots,\la_n$, where
$\la_1\leqslant\dots\leqslant \la_n$ and $\la_1+\dots+\la_n=N$.
The centralizer $\agot=\g^e$ is a Lie algebra
with the basis elements $E_{ij}^{(r)}$, where the range of indices
is defined by the inequalities $1\leqslant i,j\leqslant n$
and $\la_j-\min(\la_i,\la_j)\leqslant r<\la_j$,
with the commutation relations
\ben
\big[E^{(r)}_{ij},E^{(s)}_{kl}\big]=\de_{kj}\ts E^{(r+s)}_{i\tss l}-\de_{i\tss l}\ts E^{(r+s)}_{kj},
\een
assuming that $E^{(r)}_{ij}=0$ for $r\geqslant \la_j$.
The formulas expressing $E^{(r)}_{ij}$ in terms of the basis elements of $\gl_N$
can be found e.g. in \cite{bb:ei} and \cite{ppy:si}. It is clear from
the relations that
in the particular case $\la_1=\dots=\la_n=p$
the Lie algebra $\agot$ is isomorphic to the truncated polynomial current algebra
$\gl_n[x]/(x^p=0)$ (also known as the {\em Takiff algebra}).

Equip the Lie algebra $\agot$ with
the invariant symmetric bilinear form $(\ts\ts|\ts\ts)$ defined by
\beql{formtr}
(E_{ij}^{(0)}\tss|\tss E_{ji}^{(0)})=\la_{i},
\eeq
whereas all remaining values of the form on the basis vectors are zero.
It is understood that if $i\ne j$ then relation \eqref{formtr}
occurs only if $\la_i=\la_j$.

Consider the {\em differential algebra} $\Vc(\agot)$
which is defined as the algebra
of differential polynomials
in the variables $E^{(r)}_{ij}[s]$, where $s=0,1,2,\dots$,
while $E^{(r)}_{ij}$ ranges over the basis elements of $\agot$,
equipped with the derivation $\di$ defined by $\di\tss(E^{(r)}_{ij}[s])=E^{(r)}_{ij}[s+1]$.
We will regard $\agot$ as a subspace of $\Vc(\agot)$ by using the embedding $X\mapsto X[0]$
for $X\in\agot$.

Introduce the $\lam$-{\em bracket} on $\Vc(\agot)$ as a linear map
\ben
\Vc(\agot)\ot\Vc(\agot)\to\CC[\lam]\ot\Vc(\agot),\qquad a\ot b\mapsto \{a_{\ls}b\}.
\een
By definition, it is given by
\ben
\{X_{\ls}Y\}=[X,Y]+(X|Y)\tss\lam\qquad\text{for}\quad X,Y\in\agot,
\een
and extended to $\Vc(\agot)$ by sesquilinearity $(a,b\in\Vc(\agot))$:
\ben
\{\di\tss a_{\ls}b\}=-\ls\tss\{a_{\ls}b\},\qquad
\{a_{\ls}\di\tss b\}=(\lam+\di)\tss\{a_{\ls}b\},
\een
skewsymmetry
\ben
\{a_{\ls}b\}=-\{b_{\ts-\lam-\di\ts}a\},
\een
and the Leibniz rule $(a,b,c\in\Vc(\agot))$:
\ben
\{a_{\ls}b\tss c\}=\{a_{\ls}b\}\tss c+\{a_{\ls}c\}\tss b.
\een
The $\lam$-bracket defines the {\it affine Poisson vertex algebra\/} structure on
$\Vc$; see \cite{dskv:cw}.

Consider the following triangular decomposition of the Lie algebra $\agot$,
\beql{triang}
\agot=\n_-\oplus\h\oplus \n_+,
\eeq
where the subalgebras are defined by
\ben
\n_-=\text{span of\ }\{E^{(r)}_{ij}\ts|\ts i>j\},\quad
\n_+=\text{span of\ }\{E^{(r)}_{ij}\ts|\ts i<j\}\fand \h=\text{span of\ }\{E^{(r)}_{ii}\},
\een
with the superscript $r$ ranging over all admissible values.
Set $\p=\n_-\oplus\h$ and regard $\Vc(\p)$ as a natural differential subalgebra of
$\Vc(\agot)$. Define the differential algebra homomorphism
\ben
\rho:\Vc(\agot)\to \Vc(\p)
\een
as the identity map on all elements of $\p\subset\agot$, while for the basis elements
$E^{(r)}_{ij}\in\n_+$
we set
\beql{rhopro}
\rho\tss(E^{(r)}_{ij})=\begin{cases}1\qquad&\text{if}\quad j=i+1\fand r=\la_{i+1}-1,\\
                                     0\qquad&\text{otherwise}.
                         \end{cases}
\eeq

\bde\label{def:clawa}
The {\em classical $\Wc$-algebra $\Wc(\agot)$} associated with
$\agot$
is defined by
\ben
\Wc(\agot)=\{P\in\Vc(\p)\ |\ \rho\tss\{X_{\ls}P\}=0\quad\text{for all}\quad
X\in\n_+\}.
\een
\ede

In the case $e=0$ this agrees with the standard terminology; cf. \cite{dskv:cw}.

\bre\label{rem:pola}
One can define a more general family of classical $\Wc$-algebras $\Wc_A(\agot)$
by altering the map $\rho$ and making it dependent on
a chosen set of polynomials $A_1(u),\dots,A_{n-1}(u)$ of the form
\ben
A_i(u)=A_i^{(\la_{i+1}-\la_i)}\tss u^{\la_{i+1}-\la_i}
+\dots+A_i^{(\la_{i+1}-1)}\tss u^{\la_{i+1}-1}
\een
for $i=1,\dots,n-1$. Instead of \eqref{rhopro} one then takes
\ben
\rho\tss(E^{(r)}_{ij})=\begin{cases}A_i^{(r)}\qquad&\text{if}\quad j=i+1,\\
                                     0\qquad&\text{otherwise}.
                         \end{cases}
\een
However, it follows from the arguments below (see the proof of Theorem~\ref{thm:glncent})
that all $\Wc$-algebras
$\Wc_A(\agot)$ are isomorphic to each other provided that the leading coefficient
of each polynomial $A_i(u)$ is nonzero. For that reason we will stick to the
particular choice of such polynomials $A_i(u)=u^{\la_{i+1}-1}$ as given
in \eqref{rhopro}.
\qed
\ere

The following properties are verified by using the same arguments as in the case $e=0$; see
\cite[Sec.~3.2]{dskv:cw}.

\ble\label{lem:repcla}\quad {\rm (i)}
The mapping
\ben
(a,g)\mapsto \rho\tss\{a_{\ls}g\},\qquad a\in \n_+,\quad g\in \Vc(\p)
\een
defines a representation of the Lie conformal algebra $\CC[\di]\tss\n_+$
on $\Vc(\p)$.
\par
{\rm (ii)} The Lie conformal algebra $\CC[\di]\tss\n_+$ acts on $\Vc(\p)$
by conformal derivations so that
\ben
\rho\tss\{a_{\ls}gh\}=\rho\tss\{a_{\ls}g\}\tss h+\rho\tss\{a_{\ls}h\}\tss g
\een
for all $a\in\n_+$ and $g,h\in \Vc(\p)$.
\qed
\ele

\bpr\label{prop:pva}
The subspace
$\Wc(\agot)\subset \Vc(\p)$ is a differential subalgebra. Moreover,
$\Wc(\agot)$ is a Poisson vertex algebra
equipped with the $\lam$-bracket
\ben
(a,b)\mapsto\rho\tss\{a_{\ls}b\},\qquad a,b\in \Wc(\agot).
\vspace{-0.7cm}
\een
\qed
\epr

\section{Generators of $\Wc(\agot)$}
\label{sec:gen}

For an $n\times n$ matrix $A=[a_{ij}]$
with entries in a ring we will consider its {\em column-determinant}
defined by
\ben
\cdet\ts A=\sum_{\si\in\Sym_n}\sgn\si\cdot a_{\si(1)\ts 1}\dots a_{\si(n)\ts n}.
\een

Recall that the $\la_i$ form an increasing sequence and
combine the basis elements of the Lie algebra $\agot$ into
polynomials in a variable $u$ by
\ben
E_{ij}(u)=\begin{cases}E^{(0)}_{ij}+\dots+E^{(\la_j-1)}_{ij}\ts u^{\la_j-1}
&\text{if}\quad i\geqslant j,\\[0.4em]
E^{(\la_j-\la_i)}_{ij}\tss u^{\la_j-\la_i}+\dots+E^{(\la_j-1)}_{ij}\ts u^{\la_j-1}
&\text{if}\quad i< j.
\end{cases}
\een
For a variable $x$ the column-determinant
\ben
D_n=\cdet\ts
\begin{bmatrix}x+\la_1\tss\di+E_{11}(u)&u^{\la_2-1}&0&\dots&0\\[0.4em]
                 E_{21}(u)&x+\la_2\tss\di+E_{22}(u)&u^{\la_3-1}&\dots&0\\[0.4em]
                 \dots&\dots&\dots&\dots& \dots \\
                             \dots&\dots&\dots&\dots&u^{\la_n-1}\\[0.4em]
                             E_{n1}(u)&E_{n2}(u)&\dots&\dots&x+\la_n\tss\di+E_{nn}(u)
                \end{bmatrix}
\een
is a polynomial in $x$ of the form
\ben
x^{\tss n}+w_1(u)\tss x^{\tss n-1}+\dots+w_n(u),\qquad w_k(u)=\sum_r w^{(r)}_k\ts u^r,
\qquad w^{(r)}_k\in\Vc(\p).
\een

The following theorem is our main result. Its particular case for the element $e=0$
(that is, with $\la_1=\dots=\la_n=1$) is \cite[Thm~3.2]{mr:cw}.

\bth\label{thm:glncent}
All elements $w^{(r)}_k$ with $k=1,\dots,n$ and
\beql{conda}
\la_{n-k+2}+\dots+\la_n< r+k\leqslant\la_{n-k+1}+\dots+\la_n
\eeq
belong to the classical $\Wc$-algebra
$\Wc(\agot)$. Moreover, under these conditions
the elements $\di^{\tss s} w^{(r)}_k$
with $s=0,1,\dots$ are algebraically independent
and generate the algebra $\Wc(\agot)$.
\eth

\bex\label{ex:n2}
For $n=2$ we have
\ben
w_1(u)=E_{11}(u)+E_{22}(u),\qquad w_2(u)=E_{11}(u) E_{22}(u)-u^{\la_2-1}E_{21}(u)+\la_1\tss \di E_{22}(u),
\een
so that
\ben
w_1^{(r)}=E_{11}^{(r)}+E_{22}^{(r)},\qquad r=0,1,\dots,\la_2-1,
\een
and
\ben
w_2^{(r)}=\sum_{a+b=r} E_{11}^{(a)}E_{22}^{(b)}-E_{21}^{(r-\la_2+1)}+\la_1\tss E_{22}^{(r)}[1]
\een
for $r=\la_2-1,\dots,\la_1+\la_2-2$.
\qed
\eex

\bpf[Proof of Theorem~\ref{thm:glncent}]
Due to Lemma~\ref{lem:repcla}\ts(i), the first part of the theorem will follow if
we show that for $i=1,\dots,n-1$ the relations
\ben
\rho\{{E_{i\ts i+1}^{(t)}}_{\ls} w_k^{(r)}\}=0
\een
hold for all $t=\la_{i+1}-\la_i,\dots,\la_{i+1}-1$, assuming the conditions
on the values of $k$ and $r$ as stated in the theorem.

Set $d_i=x+\la_i\tss\di$ for $i=1,\dots,n$ and note
the commutation relations
\beql{defdx}
\big[d_i,X[s]\tss\big]=\la_i\ts X[s+1],\qquad X\in \p.
\eeq
For any $X\in \n_+$
and $Y\in \Vc(\p)$ we have the relations
\ben
\{X_{\ls}d_iY\}=(d_i+\lam\la_i)\{X_{\ls}Y\}.
\een
We will set $d_i^+=d_i+\lam\la_i$ for brevity.
For all $k\geqslant l$ we have
\ben
\{{E_{i\ts i+1}^{(t)}}_{\ls} E_{k\tss l}(u)\}=\de_{k\ts i+1}\ts E^{\tss t}_{i\tss l}(u)
-\de_{i\tss l}\ts E^{\tss t}_{k\tss i+1}(u)+\de_{t\tss 0}\ts\de_{k\ts i+1}\ts \de_{i\ts l}\ts \lam\la_i,
\een
where we set
\ben
E^{\tss t}_{ij}(u)=E^{(t)}_{ij}+\dots+E^{(\la_j-1)}_{ij}\ts u^{\la_j-t-1}.
\een
By using Lemma~\ref{lem:repcla}\ts(ii) and taking the $\lam$-bracket of
$E_{i\ts i+1}^{(t)}$ with the elements of consecutive columns of
the column-determinant $D_n$, we get
\ben
\rho\{{E_{i\ts i+1}^{(t)}}_{\ls} D_n\}
=\sum_{a=1}^{i+1} D_{n\tss a},
\een
where $D_{n\tss a}$ is the column-determinant obtained from $D_n$ as follows.
The elements $d_c$ should be replaced by the rule
$d_c\mapsto d^+_c$ for $c=1,\dots,a-1$ and column $a$ in $D_n$ should be replaced by
a new column which we now describe.

For $a=1,\dots,i-1$ we have
\ben
D_{n\tss a}=\cdet\ts
\begin{bmatrix}d^+_1+E_{11}(u)&u^{\la_2-1}&\dots&0&\dots&0\\[0.4em]
                 E_{21}(u)&d^+_2+E_{22}(u)&\dots&0&\dots&0\\[0.4em]
                 E_{31}(u)&E_{32}(u)&\dots&\vdots&\dots&0\\[0.4em]
                 \vdots&\vdots&\dots&E^{\tss t}_{i\tss a}(u)&\dots& \vdots \\
                             \vdots&\vdots&\dots&\vdots&\dots&u^{\la_n-1}\\[0.4em]
                             E_{n1}(u)&E_{n2}(u)&\dots&0&\dots&d_n+E_{nn}(u)
                \end{bmatrix},
\een
where the only nonzero entry $E^{\tss t}_{i\tss a}(u)$ in column $a$ occurs in row $i+1$.

The column-determinant $D_{n\tss i+1}$ has the same form, except that
the only nonzero entry in column $i+1$ occurring in row $i+1$ equals $u^{\la_{i+1}-t-1}$.

In the remaining column-determinant $D_{n\tss i}$ we only display
columns $i$ and $i+1$ which have the form
\ben
\cdet\ts
\begin{bmatrix}\quad\dots&0&0&\dots\quad\\
                 \quad\dots&\vdots&\vdots&\dots\quad\\
                 \quad\dots&-u^{\la_{i+1}-t-1}&u^{\la_{i+1}-1}&\dots\quad\\
                 \quad\dots&E^{\tss t}_{ii}(u)-E^{\tss t}_{i+1\ts i+1}(u)+\de_{t\tss 0}\ts\lam\tss\la_i
                             &d_{i+1}+E_{i+1\ts i+1}(u)&\dots\quad\\[0.3em]
                 \quad\dots&-E^{\tss t}_{i+2\ts i+1}(u)&E_{i+2\ts i+1}(u)&\dots\quad\\
                 \quad\dots&\vdots&\vdots&\dots\quad\\
                 \quad\dots&-E^{\tss t}_{n\ts i+1}(u)&E_{n\ts i+1}(u)&\dots\quad
                \end{bmatrix}.
\een
By expanding $D_{n\tss i}$ along columns $i$ and $i+1$ simultaneously, we find that a nonzero
contribution from this expansion can come only from the $2\times 2$ column-minors with entries
in rows $i$ and $b$ with $b\geqslant i+1$. For $b>i+1$ such a minor equals
\beql{mintwotwo}
-u^{\la_{i+1}-t-1}\tss E_{b\ts i+1}(u)+u^{\la_{i+1}-1}\tss E^{\tss t}_{b\ts i+1}(u)
\eeq
which is zero for $t=0$ and is equal to
a certain polynomial $C_b(u)$ in $u$ of degree $\la_{i+1}-2$ for $t>0$. In the latter case
the contribution to $D_{n\tss i}$ arising from this $2\times 2$ column-minor equals
an alternating sum of products of the form
\begin{multline}
\non
\big(\de_{\si(1)\ts 1}\ts d^+_1+E_{\si(1)\ts 1}(u)\big)\dots
\big(\de_{\si(i-1)\ts i-1}\ts d^+_{i-1}+E_{\si(i-1)\ts i-1}(u)\big)\ts C_b(u)\\[0.4em]
{}\times\big(\de_{\si(i+2)\ts i+2}\ts d_{i+2}+E_{\si(i+2)\ts i+2}(u)\big)\dots
\big(\de_{\si(n)\ts n}\ts d_{n}+E_{\si(n)\ts n}(u)\big)
\end{multline}
with $\si\in\Sym_n$, where each occurrence of $E_{j\ts j+1}(u)$ amongst
$E_{\si(l)\ts l}(u)$
should be
replaced with $u^{\la_{j+1}-1}$, and we assume that $E_{j\tss l}(u)=0$ for $l>j+1$.
Since the $u$-degree of $E_{j\tss l}(u)$ equals $\la_{l}-1$,
the resulting expression will be a polynomial in $x$ such that
the coefficient of $x^{\tss n-k}$ is a polynomial in $u$ whose degree
does not exceed
\ben
(\la_{n-k+1}-1)+\dots+(\la_{i-1}-1)+(\la_{i+1}-2)+(\la_{i+2}-1)+\dots+(\la_{n}-1)
\een
if $n-k\leqslant i-1$, and does not exceed
\ben
(\la_{i+1}-2)+(\la_{n-k+3}-1)+\dots+(\la_{n}-1)
\een
if $n-k\geqslant i$. In both bases the degree in $u$ is less than
\beql{sumnk}
(\la_{n-k+2}-1)+\dots+(\la_{n}-1)
\eeq
and so the contribution of such $2\times 2$ column-minors to the
expression $\rho\{{E_{i\ts i+1}^{(t)}}_{\ls} w_k^{(r)}\}$ is equal to zero
because by the conditions of the theorem the minimal
possible value of $r$ is given by \eqref{sumnk}.

Furthermore, we have
\ben
D_{n\tss i+1}=D^+_{1,\dots,i}\tss
u^{\la_{i+1}-t-1}\tss D_{i+2,\dots,n},
\een
where for a subset $K\subset \{1,\dots,n\}$ we denote by $D_K$ the principal
column-minor of $D_n$ corresponding to the rows and columns labelled by
the elements of $K$, while $D^+_K$ is obtained from $D_K$ by
replacing all elements $d_a$ with $d_a^+$.
Expand $D^+_{1,\dots,i}$
along the last row to get
\ben
D^+_{1,\dots,i}=\sum_{a=1}^{i-1} (-1)^{a+i}\tss
D^+_{1,\dots,a-1}\ts u^{\la_{a+1}-1}\dots u^{\la_{i}-1}
\ts E_{i\tss a}(u)+D^+_{1,\dots,i-1}\tss \big(d^+_{i}+E_{ii}(u)\big).
\een
This shows that if $t>0$ then $\rho\{{E_{i\ts i+1}^{(t)}}_{\ls} D_n\}$
can be written as the expression
\begin{multline}
\non
\sum_{a=1}^{i} (-1)^{a+i}\tss D^+_{1,\dots,a-1}\ts u^{\la_{a+1}-1}\dots u^{\la_{i}-1}
\big(u^{\la_{i+1}-t-1}\tss E_{i\tss a}(u)-u^{\la_{i+1}-1}\tss E^{\tss t}_{i\tss a}(u)\big)\tss
\tss D_{i+2,\dots,n}\\
{}-D^+_{1,\dots,i-1}\ts\Big(
\big(u^{\la_{i+1}-t-1}\tss E_{i+1\tss i+1}(u)-u^{\la_{i+1}-1}\tss E^{\tss t}_{i+1\tss i+1}(u)\big)
-u^{\la_{i+1}-t-1}\tss (d^+_i-d^{}_{i+1})\tss \Big)\\
{}\times D_{i+2,\dots,n},
\end{multline}
while for the value $t=0$ occurring in the case $\la_i=\la_{i+1}$ we have
\ben
\rho\{{E_{i\ts i+1}^{(0)}}_{\ls} D_n\}=D^+_{1,\dots,i-1}\ts u^{\la_{i+1}-1}\ts
(d_i-d_{i+1})\ts D_{i+2,\dots,n}.
\een
Similar to the argument above with elements \eqref{mintwotwo}, we can see that
the expression $\rho\{{E_{i\ts i+1}^{(t)}}_{\ls} D_n\}$
is a polynomial in $x$ such that
the coefficient of $x^{\tss n-k}$ is a polynomial in $u$ whose degree
is less than the number in \eqref{sumnk}. This completes the proof
of the first part of the theorem.

As a next step, we will show that the elements $\di^{\tss s} w^{(r)}_k$
are algebraically independent. We will adapt the corresponding argument
used in the proof of \cite[Theorem~3.14]{dskv:cw} and introduce
the differential polynomial degree on $\Vc(\p)$ by setting the degree of $X[s]$
to be equal to $s+1$ for any nonzero $X\in\p$. The minimal degree components
of the elements $w^{(r)}_k$ with $r=q+\la_{n-k+2}+\dots+\la_n-k+1$ are given by the formulas
\ben
w^{(r)}_k=(-1)^{k-1}\ts\sum_{p=k}^n E_{p\ts p-k+1}^{(q+\la_n-\la_{n-k+1}+\dots+\la_{p+1}-\la_{p-k+2})}
+\text{\ higher degree terms}
\een
for $q=0,1,\dots,\la_{n-k+1}-1$. By applying $s$ times the derivation $\di$ to both sides
we get the respective minimal degree components
of the elements $\di^{\tss s} w^{(r)}_k$. Their algebraic independence now follows
from the observation that all elements of the form
\ben
\sum_{p=k}^n E_{p\ts p-k+1}^{(q+\la_n-\la_{n-k+1}+\dots+\la_{p+1}-\la_{p-k+2})}[s]
\een
with $s=0,1,\dots$ are algebraically independent.

Finally, we will show that the algebra $\Wc(\agot)$ is generated by the elements
$\di^{\tss s} w^{(r)}_k$. Expanding $D_n$
along the last row we get the expression
\ben
D_n=D_{1,\dots,n-1}\tss \big(d_{n}+E_{nn}(u)\big)
+\sum_{k=2}^{n} (-1)^{k-1}\tss
D_{1,\dots,n-k}\ts u^{\la_{n-k+2}-1}\dots u^{\la_{n}-1}
\ts E_{n\tss n-k+1}(u).
\een
Hence
for $k=1,\dots,n$ and $r=q+\la_{n-k+2}+\dots+\la_n-k+1$ we can write
\beql{wreq}
w^{(r)}_k=(-1)^{k-1}\ts E^{(q)}_{n\tss n-k+1}+R,
\qquad q=0,1,\dots,\la_{n-k+1}-1,
\eeq
where $R$ is a polynomial in the variables of the form
$E^{(p)}_{n\tss m}[s]$ with $m=n-k+2,\dots,n$ and the variables
$E^{(p)}_{l\tss m}[s]$ with $n>l\geqslant m\geqslant 1$.

Now suppose that an element $P\in\Vc(\p)$ belongs to the subalgebra $\Wc(\agot)\subset\Vc(\p)$.
Applying relations \eqref{wreq} consecutively with $k=n,n-1,\dots,1$
we can write $P$ as a polynomial in the new variables
$\di^{\tss s} w^{(r)}_k$ with $k=1,\dots,n$ and $s\geqslant 0$ satisfying
conditions \eqref{conda}, together with $E^{(p)}_{l\tss m}[s]$ with $n>l\geqslant m\geqslant 1$.
Take $X=E_{i\ts n}^{(\la_n-t-1)}$ for $i\in\{1,\dots,n-1\}$
and $t\in\{0,1,\dots,\la_i-1\}$
in Definition~\ref{def:clawa} and observe that the images of the $\lam$-brackets
of $X$ with the new variables under the map $\rho$ are equal to zero, except for
the variables $E_{n-1\ts i}^{(t)}[s]$. Since
\ben
\{{E_{i\ts n}^{(\la_n-t-1)}}_{\ls} E_{n-1\ts i}^{(t)}[s]\}=-(\lam+\di)^s \ts E_{n-1\ts n}^{(\la_n-1)},
\een
by applying Lemma~\ref{lem:repcla}\ts(ii) we find
that the property $\rho\tss\{X_{\ls}P\}=0$ implies the relation
\ben
\sum_{s\geqslant 0} \lam^s\ts\frac{\di P}{\di E_{n-1\ts i}^{(t)}[s]}=0
\een
which means that $P$ does not depend on the variables $E_{n-1\ts i}^{(t)}[s]$
with $i=1,\dots,n-1$. Repeating the same argument for the elements $X\in\n_+$
of the form $X=E_{i\ts j}^{(\la_j-t-1)}$ with $j=n-1,\dots,2$ and $i\in\{1,\dots,j-1\}$,
we may conclude that $P$ does not depend on the variables
$E^{(p)}_{l\tss m}[s]$ with $n>l\geqslant m\geqslant 1$ so that
$P$ is a polynomial in the variables $\di^{\tss s} w^{(r)}_k$, as required.
\epf

\section{Miura map and center at the critical level}
\label{sec:crit}

Here we will show that the classical $\Wc$-algebra $\Wc(\agot)$ is isomorphic to
the center at the critical level
of the affine vertex algebra $V(\agot)$; cf.~\cite{ff:ak},
\cite[Ch.~4]{f:lc}. We will do this by relying on the work \cite{m:cc} and
providing an explicit correspondence
between the generators of both algebras. We will use a Miura-type map
on the $\Wc$-algebra side and a Harish-Chandra-type isomorphism on the vertex algebra side.

Denote by $\Vc(\h)$ the differential subalgebra of $\Vc(\p)$, generated by the
elements $E_{i\tss i}^{(r)}$ with $i=1,\dots,n$ and $r=0,\dots,\la_i-1$.
Let
\ben
\vp:\Vc(\p)\to\Vc(\h)
\een
denote the homomorphism of differential algebras defined on the generators
as the projection $\p\to\h$ with the kernel $\n_-$.

\bpr\label{prop:chev}
The restriction of the homomorphism $\vp$ to the classical $\Wc$-algebra
$\Wc(\agot)$ is injective.
\epr

\bpf
It will be sufficient to verify that the images of all generators
$\di^{\tss s} w^{(r)}_k\in\Wc(\agot)$ with $k=1,\dots,n$ and $s\geqslant 0$ satisfying
conditions \eqref{conda} are algebraically independent in $\Vc(\h)$.
The images $\overline w^{\tss(r)}_k=\vp(w^{(r)}_k)$ are found by writing the product
\ben
\big(x+\la_1\tss\di+E_{1\tss 1}(u)\big)\dots \big(x+\la_n\tss\di+E_{n\tss n}(u)\big)
\een
as a polynomial in $x$ of the form
\ben
x^{\tss n}+\overline w_1(u)\tss x^{\tss n-1}+\dots+\overline w_n(u),\qquad
\overline w_k(u)=\sum_r \overline w^{\tss(r)}_k\ts u^r.
\een
Introduce a grading on $\Vc(\h)$ by setting the degree of $X[s]$
to be equal to $s$ for any nonzero $X\in\h$. The proposition will follow if we show that
the minimal degree components
$\di^{\tss s} v^{(r)}_k$ of the respective
elements $\di^{\tss s}\tss\overline w^{\tss(r)}_k$ are algebraically independent.
These components are found by the formulas
\beql{polv}
\di^{\tss s} v^{(r)}_k=\sum_{i_1<\dots<i_k}\ts\sum_{r_1+\dots+r_k=r}\ts\sum_{s_1+\dots+s_k=s}
\frac{s!}{s_1!\dots s_k!}\ts E_{i_1i_1}^{(r_1)}[s_1]\dots E_{i_ki_k}^{(r_k)}[s_k].
\eeq
We will verify that the differentials of these polynomials are linearly independent.
Introduce particular orderings on the set of the variables $E_{ii}^{(r)}[s]$
and on the set of polynomials $\di^{\tss s} v^{(r)}_k$ as follows. First, if $s<s'$
then we set
$E_{i\tss i}^{(r)}[s]\prec E_{j\tss j}^{(p)}[s']$ and
$\di^{\tss s} v^{(r)}_k\prec \di^{\tss s'} v^{(p)}_l$
for all admissible values of the remaining parameters. Formulas \eqref{polv} imply
that if the orderings within the sets of elements with a given $s$
are chosen consistently for different values of $s$, then
the Jacobian matrix will have a block-diagonal form with identical diagonal blocks.
Therefore, it will be sufficient to consider the Jacobian matrix
corresponding to the elements with $s=0$. In this case
formula \eqref{polv} gives
\beql{polvzero}
v^{(r)}_k=\sum_{i_1<\dots<i_k}\ts\sum_{r_1+\dots+r_k=r}\ts
E_{i_1i_1}^{(r_1)}\dots E_{i_ki_k}^{(r_k)}.
\eeq
Now list the variables in a particular order:
\ben
E_{n\tss n}^{(\la_n-1)},\dots,E_{11}^{(\la_1-1)},E_{n\tss n}^{(\la_n-2)},\dots,E_{11}^{(\la_1-2)},
\dots,E_{n\tss n}^{(0)},\dots,E_{11}^{(\la_1-\la_n)},
\een
assuming that any variables with negative superscripts are excluded from the list.
This means, in particular, that the last segment of the sequence
has the form $E_{n\tss n}^{(0)},\dots,E_{l+1\tss l+1}^{(0)}$, if for certain
index $l$ we have $\la_l<\la_{l+1}=\dots=\la_n$.
Similarly, list the polynomials by
\ben
v^{(\la_n-1)}_1,\dots,v^{(\la_n+\dots+\la_1-n)}_n,
v^{(\la_n-2)}_1,\dots,v^{(\la_n+\dots+\la_1-n-1)}_n,\dots,
v^{(0)}_1,\dots,v^{(\la_{n-1}+\dots+\la_{1}-n+1)}_n,
\een
where a polynomial $v^{(r)}_k$ is excluded from the list if conditions
\eqref{conda} do not hold. Note that the corresponding Jacobian matrix is square
of size $N=\la_1+\dots+\la_n$. Introduce multiplicities $q_1,\dots,q_l$
of the decreasing sequence $(\la_n,\la_{n-1},\dots,\la_1)$ by
\ben
\la_n=\dots=\la_{n-q_1+1}>\la_{n-q_1}=\dots=\la_{n-q_1-q_2+1}>\dots>
\la_{q_l}=\dots=\la_1,
\een
so that the $q_i$ are positive integers with $q_1+\dots+q_l=n$. First consider
the upper left block of the Jacobian matrix of size $n\times n$. By
formulas \eqref{polvzero}, the elements
\ben
v^{(\la_n-1)}_1,\dots,v^{(\la_n+\dots+\la_{n-q_1+1}-q_1)}_{q_1}
\een
are the respective elementary symmetric polynomials in the variables
\ben
E_{n\tss n}^{(\la_n-1)},\dots,E_{n-q_1+1\ts n-q_1+1}^{(\la_{n-q_1+1}-1)}.
\een
It is well-known that the corresponding $q_1\times q_1$ Jacobian matrix is non-degenerate
which is equivalent to the algebraic independence of the
elementary symmetric polynomials. Similarly, each of the next $q_2$ elements
\ben
v^{(\la_n+\dots+\la_{n-q_1}-q_1-1)}_{q_1+1},\dots,v^{(\la_n+\dots+\la_{n-q_1-q_2+1}-q_1-q_2)}_{q_1+q_2}
\een
is the respective elementary symmetric polynomial in the variables
\ben
E_{n-q_1\ts n-q_1}^{(\la_{n-q_1}-1)},\dots,E_{n-q_1-q_2+1\ts n-q_1-q_2+1}^{(\la_{n-q_1-q_2+1}-1)}
\een
multiplied by the product $E_{n\tss n}^{(\la_n-1)}\dots E_{n-q_1+1\ts n-q_1+1}^{(\la_{n-q_1+1}-1)}$.
Therefore, the next $q_2\times q_2$ diagonal block
of the Jacobian matrix is also non-degenerate. A similar structure is clearly retained
by all remaining $q_i\times q_i$ diagonal blocks. In particular, each of the last $q_l$ elements
\ben
v^{(\la_n+\dots+\la_{q_l}-n+q_l-1)}_{n-q_l+1},\dots,v^{(\la_n+\dots+\la_1-n)}_{n}
\een
is the respective elementary symmetric polynomial in the variables
$
E_{q_l q_l}^{(\la_{q_l}-1)},\dots,E_{1\tss 1}^{(\la_{1}-1)}
$
multiplied by the product of the first $n-q_l$ variables
$E_{n\tss n}^{(\la_n-1)}\dots E_{q_l+1\ts q_l+1}^{(\la_{q_l+1}-1)}$.

By formulas \eqref{polvzero},
each of the next $\la_1-1$ diagonal $n\times n$ blocks of the Jacobian matrix
exhibits a similar structure. Namely, for each
$i=1,\dots,l$ their $q_i\times q_i$
diagonal sub-blocks coincide with
the respective $q_i\times q_i$ sub-block in the first block of size $n\times n$.

As a next step, consider the submatrix of the Jacobian matrix
of size $(n-q_l)\times (n-q_l)$
corresponding to the subsequence of variables
\ben
E_{n\tss n}^{(\la_n-\la_1-1)},\dots,E_{q_l+1\ts q_l+1}^{(\la_{q_l+1}-\la_1-1)}
\een
and the polynomials
\ben
v^{(\la_n-\la_1-1)}_1,\dots,v^{(\la_n+\dots+\la_{q_l+1}-\la_1-n+q_l)}_{n-q_l}.
\een
It follows from \eqref{polvzero} that
this submatrix has a block-triangular form with the diagonal blocks of sizes
$q_i\times q_i$ for $i=1,\dots,l-1$ which coincide with the respective
sub-blocks of the $n\times n$ submatrix considered above. The argument
continues in the same way for the remaining variables and polynomials, and the same observation
applies to the remaining part of the Jacobian matrix. Since
all diagonal blocks are non-degenerate we may conclude that
so is the full matrix.
\epf

\bre\label{rem:algin}(i)\ \
The above argument provides another proof
of the algebraic independence of the generators
$\di^{\tss s} w^{(r)}_k$ of the algebra $\Wc(\agot)$ given in Theorem~\ref{thm:glncent}.
\par
(ii)\ \  Proposition~\ref{prop:chev} allows one to regard the classical $\Wc$-algebra
$\Wc(\agot)$ as a differential subalgebra of $\Vc(\h)$ generated by
the elements $\overline w^{\tss(r)}_k$ with $k=1,\dots,n$ and $r$ satisfying
conditions \eqref{conda}. It would be interesting to find an intrinsic characterization of this
subalgebra. This could possibly involve some
version of {\em screening operators} like in the case $e=0$; see \cite[Ch.~8]{f:lc}.
In this case the embedding $\Wc(\agot)\hra\Vc(\h)$
is the {\em Miura map}, {\it  loc. cit.}
\qed
\ere

Now we recall the definition of the center at the critical level
associated with the centralizer $\agot$ from \cite{ap:qm}
and a construction of its generators \cite{m:cc}.
The affine Kac--Moody algebra $\wh\agot$
is the central
extension
\ben
\wh\agot=\agot\tss[t,t^{-1}]\oplus\CC \vac,
\een
where $\agot[t,t^{-1}]$ is the Lie algebra of Laurent
polynomials in $t$ with coefficients in $\agot$. For any $r\in\ZZ$ and $X\in\g$
we will write $X[r]=X\ts t^r$. The commutation relations of the Lie algebra $\wh\agot$
have the form
\ben
\big[X[r],Y[s]\big]=[X,Y][r+s]+r\ts\de_{r,-s}\langle X,Y\rangle\ts \vac,
\qquad X, Y\in\agot,
\een
and the element $\vac$ is central in $\wh\agot$.
Here the invariant symmetric bilinear form $\langle\ts\ts,\ts\rangle$ is
different from \eqref{formtr} and
defined by the formulas
\ben
\big\langle E_{ii}^{(0)},E_{jj}^{(0)}\big\rangle=
\min(\la_i,\la_j)-\de_{ij}\big(\la_1+\dots+\la_{i-1}+(n-i+1)\la_i\big),
\een
and if $\la_i=\la_j$ for some $i\ne j$ then
\ben
\big\langle E_{ij}^{(0)},E_{ji}^{(0)}\big\rangle=-\big(\la_1+\dots+\la_{i-1}+(n-i+1)\la_i\big),
\een
whereas all remaining values of the form on the basis vectors are zero.
The {\em vacuum module at the critical level}
over $\wh\agot$
is the quotient
\ben
V(\agot)=\U(\wh\agot)/\Ir,
\een
where $\Ir$ is the left ideal of $\U(\wh\agot)$ generated by $\agot[t]$
and the element $\vac-1$. The vacuum module
is a vertex algebra and its {\em center} is defined as the subspace
\ben
\z(\wh\agot)=\{v\in V(\agot)\ |\ \agot[t]\tss v=0\}.
\een
The center is a commutative associative algebra which can be regarded as
a subalgebra of $\U\big(t^{-1}\agot[t^{-1}]\big)$.
This subalgebra is invariant with respect to the
{\em translation operator}
$T$ which is
the derivation of the algebra $\U\big(t^{-1}\agot[t^{-1}]\big)$
whose action on the generators is given by
\ben
T:X[r]\mapsto -r\tss X[r-1],\qquad X\in\agot, \quad r<0.
\een

By \cite[Thm~1.4]{ap:qm},
there exists a {\em complete set of Segal--Sugawara vectors} $S_1,\dots,S_N\in\z(\wh\agot)$,
which means that
all translations $T^r S_l$ with $r\geqslant 0$ and $l=1,\dots,N$
are algebraically independent and
any element of $\z(\wh\agot)$ can be written as a polynomial
in the shifted vectors; that is,
\ben
\z(\wh\agot)=\CC[T^{\tss r}S_l\ |\ l=1,\dots,N,\ \ r\geqslant 0].
\een
In the case $e=0$ this reduces to the Feigin--Frenkel theorem in type $A$~\cite{ff:ak, f:lc}.

To produce a complete set of Segal--Sugawara vectors,
for all $i,j\in\{1,\dots,n\}$ introduce polynomials in a variable $z$
with coefficients in this algebra by
\ben
E_{ij}(z)=\begin{cases}E^{(0)}_{ij}[-1]+\dots+E^{(\la_j-1)}_{ij}[-1]\ts z^{\la_j-1}
&\text{if}\quad i\geqslant j,\\[0.4em]
E^{(\la_j-\la_i)}_{ij}[-1]\tss z^{\la_j-\la_i}+\dots+E^{(\la_j-1)}_{ij}[-1]\ts z^{\la_j-1}
&\text{if}\quad i< j.
\end{cases}
\een
For another variable $x$
write the column-determinant
\ben
\cdet\left[\begin{matrix}
x+\la_1\tss T+E_{11}(z)&E_{12}(z)&\dots&E_{1n}(z)\\[0.2em]
E_{21}(z)&x+\la_2\tss T+E_{22}(z)&\dots&E_{2n}(z)\\
\vdots&\vdots& \ddots&\vdots     \\
E_{n1}(z)&E_{n2}(z)&\dots&x+\la_n\tss T+E_{n\tss n}(z)
                \end{matrix}\right]
\een
as a polynomial in $x$ with coefficients in $V(\agot)[z]$,
\ben
x^n+\phi_1(z)\tss x^{n-1}+\dots+\phi_n(z),\qquad \phi_k(z)=\sum_r\phi^{(r)}_k\ts z^r.
\een
By
the main result of \cite{m:cc}
the coefficients $\phi^{(r)}_k$ with $k=1,\dots,n$ and
$r$ satisfying \eqref{conda}
belong to the center $\z(\wh\agot)$
of the vertex algebra $V(\agot)$.
Moreover, they form a complete set of
Segal--Sugawara vectors for the Lie algebra $\agot$.

Denote by $\U\big(t^{-1}\agot[t^{-1}]\big)_0$ the zero weight component
of the algebra $\U\big(t^{-1}\agot[t^{-1}]\big)$ with respect to the adjoint action
of the abelian subalgebra of $\agot$ spanned by $E^{(0)}_{11},\dots,E^{(0)}_{nn}$.
Recall the triangular decomposition \eqref{triang} of $\agot$
and observe that the projection
\ben
\f:\U\big(t^{-1}\agot[t^{-1}]\big)_0\to \U\big(t^{-1}\h[t^{-1}]\big)
\een
to the first summand in the direct sum decomposition
\ben
\U\big(t^{-1}\agot[t^{-1}]\big)_0= \U\big(t^{-1}\h[t^{-1}]\big)\oplus
\Big(\U\big(t^{-1}\agot[t^{-1}]\big)_0\cap \U\big(t^{-1}\agot[t^{-1}]\big)\ts t^{-1}\n_+[t^{-1}]\Big)
\een
is an algebra homomorphism.
Note that $\z(\wh\agot)$ is a subalgebra of $\U\big(t^{-1}\agot[t^{-1}]\big)_0$.

\bpr\label{prop:hch}
The restriction of the homomorphism $\f$ to the subalgebra $\z(\wh\agot)$ is injective.
\epr

\bpf
The images of the generators of the algebra $\z(\wh\agot)$ under
the homomorphism $\f$ are readily found from their definition.
Namely, write the product
\ben
\big(x+\la_1\tss T+E_{1\tss 1}(z)\big)\dots \big(x+\la_n\tss T+E_{n\tss n}(z)\big)
\een
as a polynomial in $x$ of the form
\ben
x^{\tss n}+\overline \phi_1(z)\tss x^{\tss n-1}+\dots+\overline \phi_n(z),\qquad
\overline \phi_k(z)=\sum_r \overline \phi^{\tss(r)}_k\ts z^r.
\een
The proposition will follow if we show that the elements
$T^s\tss\overline \phi^{\tss(r)}_k\in \U\big(t^{-1}\h[t^{-1}]\big)$, where
$s\geqslant 0$ and $k=1,\dots,n$ with $r$ satisfying conditions \eqref{conda},
are algebraically independent. However, we have an isomorphism
of differential algebras
\ben
\Vc(\h)\to\U\big(t^{-1}\h[t^{-1}]\big),\qquad X[s]\mapsto s!\tss X[-s-1]
\een
for $s\geqslant 0$ and $X\in\h$, so that the derivation $\di$ corresponds to $T$.
Under this isomorphism, we have $\overline w^{\tss(r)}_k\mapsto \overline \phi^{\tss(r)}_k$.
It remains to note that the required property was already established in the proof of
Proposition~\ref{prop:chev} for the
elements $\di^{\tss s}\tss\overline w^{\tss(r)}_k$.
\epf

By the arguments used in the proof of Proposition~\ref{prop:hch} we have the following
isomorphism; cf. \cite[Theorem~8.1.5]{f:lc}.

\bco\label{cor:isomdi}
We have a differential algebra isomorphism
\ben
\z(\wh\agot)\to\Wc(\agot),\qquad \phi^{(r)}_k\mapsto w^{(r)}_k,\qquad T\mapsto \di,
\een
where $k=1,\dots,n$ and $r$ satisfies \eqref{conda}.
\qed
\eco

\section*{Acknowledgements}

E.R. wishes to thank the Sydney Mathematical Research Institute and
the School of Mathematics and Statistics at the University of Sydney for their
warm hospitality and the
SMRI International Visitor Program for the financial support.


\begin{thebibliography}{99}

\bibitem{ap:qm}
{T.~Arakawa} and {A.~Premet},
{\it Quantizing Mishchenko--Fomenko subalgebras for centralizers via affine $W$-algebras},
Trans.~Mosc.~Math.~Soc. {\bf 78} (2017), 217--234.

\bibitem{bb:ei}
{J. Brown and J. Brundan},
{\it Elementary invariants for centralizers of nilpotent matrices},
J. Aust. Math. Soc. {\bf 86} (2009), 1--15.

\bibitem{dskv:cw}
A. De Sole, V. G. Kac and D. Valeri,
{\it Classical $\Wc$-algebras and generalized Drinfeld--Sokolov bi-Hamiltonian
systems within the theory of Poisson vertex algebras},
Comm. Math. Phys. {\bf 323} (2013), 663--711.

\bibitem{dskv:cwgln}
A. De Sole, V. G. Kac and D. Valeri,
{\it Classical $\Wc$-algebras for $\gl_N$ and associated integrable Hamiltonian hierarchies},
Comm. Math. Phys. {\bf 348} (2016), 265--319.

\bibitem{dskv:ca}
A. De Sole, V. G. Kac and D. Valeri,
{\it Classical affine $\Wc$-algebras and the associated integrable
Hamiltonian hierarchies for classical Lie algebras},
Comm. Math. Phys. {\bf 360} (2018), 851--918.

\bibitem{ds:la}
{V. G. Drinfeld and V. V. Sokolov},
{\it Lie algebras and equations of Korteweg-de Vries type},
J. Sov. Math. {\bf 30} (1985), 1975--2036.

\bibitem{ff:ak}
B. Feigin and E. Frenkel,
{\it Affine Kac--Moody algebras at the critical level
and Gelfand--Dikii algebras},
Int. J. Mod. Phys. A{\bf 7}, Suppl. 1A (1992), 197--215.

\bibitem{f:lc}
E. Frenkel,
{\it Langlands correspondence for loop groups}, Cambridge
Studies in Advanced Mathematics,
103. Cambridge University Press, Cambridge, 2007.

\bibitem{m:so}
{A. Molev},
{\it Sugawara operators for classical Lie algebras}.
Mathematical Surveys and Monographs
229. AMS, Providence, RI, 2018.

\bibitem{m:cc}
{A. Molev},
{\it Center at the critical level for centralizers in type $A$},
{\tt arXiv:1904.12520}.

\bibitem{mr:cw}
A. I. Molev and E. Ragoucy,
{\it Classical W-algebras in types $A, B, C, D$ and $G$},
Comm. Math. Phys. {\bf 336} (2015), 1053--1084.

\bibitem{ppy:si}
D. Panyushev, A. Premet and O. Yakimova,
{\it On symmetric invariants of centralisers in reductive Lie algebras},
J. Algebra {\bf 313} (2007), 343--391.


\end{thebibliography}
\end{document}